\RequirePackage{fix-cm}
\documentclass[smallextended]{svjour3}       % onecolumn (second format)
\smartqed  % flush right qed marks, e.g. at end of proof
\usepackage{amsmath, amscd, amsfonts, amssymb, graphicx, color}
\begin{document}

\title{Mond and Pe$\check{c}$ari$\acute{c}$ inequality for $h$-convex functions with applications  %\thanks{Grants or other notes
%about the article that should go on the front page should be
%placed here. General acknowledgments should be placed at the end of the article.}
}
%\subtitle{Do you have a subtitle?\\ If so, write it here}

%\titlerunning{Short form of title}        % if too long for running head

\author{Ismail Nikoufar and Davuod Saeedi}

%\authorrunning{Short form of author list} % if too long for running head

\institute{I. Nikoufar (Corresponding Author) \at
              Department of Mathematics, Payame Noor University, P.O. Box 19395-3697 Tehran, Iran \\
            %  Tel.: +123-45-678910\\
             % Fax: +123-45-678910\\
              \email{nikoufar@pnu.ac.ir}           %  \\
%             \emph{Present address:} of F. Author  %  if needed
%\authorrunning{Short form of author list} % if too long for running head
\\
D. Saeedi \at
              Department of Mathematics, Payame Noor University, P.O. Box 19395-3697 Tehran, Iran \\
            %  Tel.: +123-45-678910\\
             % Fax: +123-45-678910\\
              \email{dsaeedi3961@gmail.com}           %  \\
%             \emph{Present address:} of F. Author  %  if needed
}

\date{Received: date / Accepted: date}
% The correct dates will be entered by the editor

\maketitle

\begin{abstract}
In this paper, we prove
an operator version of the Jensen's inequality and its converse for $h$-convex functions.
We provide a refinement of the Jensen type inequality for $h$-convex functions.
Moreover, we prove the Hermite-Hadamard's type inequality and a multiple operator
version of the Jensen's inequality for $h$-convex functions.
In particular, a result for convex, $P$-class, $s$-convex, Godunova-Levin, and  $s$-Godunova-Levin functions
can be deduced.

% \PACS{PACS code1 \and PACS code2 \and more}

\end{abstract}

\section{Introduction}
Throughout this paper, we denote by
$\mathcal{H}$ a Hilbert space and by $B(\mathcal{H})$ the algebra of all bounded linear operators on
$\mathcal{H}$.
The subalgebra of all self-adjoint operators in $B(\mathcal{H})$ is denoted by $B_{sa}(\mathcal{H})$.
An operator $A$ in $B_{sa}(\mathcal{H})$ is positive whenever $\langle Ax, x\rangle \geq 0$
for all $x\in\mathcal{H}$ and we write $A \geq 0$.
W denote by $Sp(A)$ the spectrum of an operator $A\in B(\mathcal{H})$.

The convexity of functions is an important issue in many fields of science, for instance in economy and optimization.
A function $f:\Bbb{I} \to \Bbb{R}$, $\Bbb{I}\subseteq \Bbb{R}$ is convex whenever the following inequality
$$
f(\lambda u + (1 - \lambda)v)\leq \lambda f(u) + (1- \lambda)f(v)
$$
holds for all $u, v \in \Bbb{I}$ and for all $\lambda\in [0,1]$ and the function $f : \Bbb{I} \to\Bbb{R}$ is concave whenever $-f$ is convex.

In 1979, Breckner \cite{Breckner} introduced the class of $s$-convex functions in the second
sense.
A function $f : [0,\infty) \to\Bbb{R}$ is $s$-convex in the second sense whenever
\begin{equation}%\label{s-convex}
f(\lambda u + (1-\lambda)v)\leq \lambda^sf(u) + (1-\lambda)^sf(v)
\end{equation}
holds for all $u, v \in [0,\infty)$, for all $\lambda\in [0, 1]$ and for some fixed $s \in (0, 1]$.
Note that all $s$-convex functions in the second sense are non-negative.
Hudzik and Maligranda (1994) \cite{Hudzik-AEqMath1994} remarked two senses of $s$-convexity of real-valued functions are known in
the literature.
A function $f : [0,\infty) \to\Bbb{R}$ is $s$-convex in the first sense if
\begin{equation}%\label{s-convex00}
f(\alpha u + \beta v)\leq \alpha^sf(u) + \beta^sf(v)
\end{equation}
holds for all $u, v \in [0,\infty)$ and for all $\alpha,\beta\geq 0$ with $\alpha^s+\beta^s=1$ and for some fixed $s \in (0, 1]$.
There is an identity between the class of 1-convex functions and the class of convex functions.
Indeed, the $s$-convexity means just the convexity when $s = 1$, no mater in the first sense or in the second sense.
For more details and examples on $s$-convex functions we refer to see \cite{Dragomir-DM1999,Dragomir-book2002,Kirmaci-AMC2007,Hudzik-AEqMath1994,Nikoufar-CAOT2021,Pinheiro-AEqMath2007,Set-CA2012}.

In 1985, Godunova and Levin (see \cite{Godunova1985} or \cite{Mitrinovic}) introduced
the class of Godunova-Levin functions.
A function $f:\Bbb{I}\to\Bbb{R}$ is a Godunova-Levin function on $\Bbb{I}$ if
\begin{equation}%\label{pi-eq}
f(\lambda u+(1-\lambda)v)\leq \frac{f(u)}{\lambda}+\frac{f(v)}{1-\lambda},
\end{equation}
where $u,v\in\Bbb{I}$ and $\lambda\in [0,1]$.
Note that all non-negative monotonic and non-negative convex functions
belong to this class \cite{Dragomir-SJM1995}.
The function $f$ is $s$-Godunova-Levin type if
\begin{equation}%\label{pi-eq}
f(\lambda u+(1-\lambda)v)\leq \frac{f(u)}{\lambda^s}+\frac{f(v)}{(1-\lambda)^s},
\end{equation}
where $u,v\in\Bbb{I}$ and $\lambda\in [0,1]$.

In 1999, Pearce and Rubinov \cite{Pearce1999} introduced a new class of convex functions
which is called $P$-class functions.
A function $f:\Bbb{I}\to\Bbb{R}$ is a $P$-class function on $\Bbb{I}$ if
\begin{equation}%\label{pi-eq00}
f(\lambda u+(1-\lambda)v)\leq f(u)+f(v),
\end{equation}
where $u,v\in\Bbb{I}$ and $\lambda\in [0,1]$.
Some properties of $P$-class functions can be found in \cite{Dragomir-SJM1995,Dragomir-BAMS1998,Nikoufar-Filomat2020}.

In 2007, in order to unify the above concepts for functions of real variable
Varo$\check{s}$anec \cite{S.Var-2007} introduced a wide class of functions
the so called $h$-convex functions which generalizes
convex, $s$-convex, Godunova-Levin, and $P$-class functions.
A non-negative function $f:\Bbb{I}\to\Bbb{R}$ is an $h$-convex function on $\Bbb{I}$ if
\begin{equation}\label{pi-eq}
f(\lambda u+(1-\lambda)v)\leq h(\lambda)f(u)+h(1-\lambda)f(v),
\end{equation}
where $h$ is a non-negative function defined on the real interval $\Bbb{J}$,
$u,v\in\Bbb{I}$ and $\lambda\in [0,1]\subseteq \Bbb{J}$.
For more results and generalizations regarding $h$-convexity, we refer the readers to see
\cite{Bombardelli2009,Dragomir-MM2015,Hazy2011,Olbrys2015}.

Exponentially convex (concave) functions are related to the convex functions
and have appeared as significant generalization of the convex functions \cite{Bernstein1929}.
%Exponentially concave functions have a significant impact in information theory \cite{Alirezaei2018} .
The concept of strongly exponentially general convex functions was introduced in \cite{Noor2020}.
These functions have some nice properties, which convex functions enjoy.

Jensen's inequality for convex  functions  is  one  of  the  most  important
result in the theory of inequalities and for appropriate choices of the function
many other famous inequalities are particular cases of this inequality.
An operator version of
the Jensen inequality for a convex function has been proved
by Mond and Pe$\check{c}$ari$\acute{c}$ as follows (\cite{Mond-HJM1993}, \cite{Furuta-Zagreb2005}):

\begin{theorem}\label{thm0}
Let $f:[m, M] \to \Bbb{R}$ be a continuous convex function.
Then,
\begin{equation}\label{eq0-thm1}
f(\langle Ax, x\rangle) \leq \langle f(A)x,x\rangle
\end{equation}
for every $x\in \mathcal{H}$ with $\langle x,x\rangle = 1$ and every self-adjoint operator $A$ such that $mI\leq A \leq MI$.
\end{theorem}

In this paper, we prove some inequalities for
self-adjoint operators on a Hilbert space including
an operator version of the Jensen's inequality and its converse.
We provide a refinement of the Jensen type inequality for $h$-convex functions.
Moreover, we prove the Hermite-Hadamard's type inequality and a multiple operator
version of the Jensen's inequality for $h$-convex functions.
In particular, a result for convex, $P$-class, $s$-convex, Godunova-Levin, and  $s$-Godunova-Levin functions
will be proved.

%%%%%%%%%%%%%%%%%%%%%%%%%%%%%%%%%%%%
%
%%%%%%%%%%%%%%%%%%%%%%%%%%%%%%%%%%%%

\section{Mond and Pe$\check{c}$ari$\acute{c}$ inequality for $h$-convex functions}
We indicate that an operator version of the Jensen inequality for $h$-convex functions still holds
similar to that Mond and Pe$\check{c}$ari$\acute{c}$ considered for convex functions.

For a non-negative and non-zero function $h:K\subset \Bbb{R}\to\Bbb{R}$, we define
$$
M_K(h)=\inf_{t\in K}\frac{h(t)}{t}
$$
and call it the Jensen's coefficient for $h$-convex functions on $K$.
We consider all functions $h$ such that the Jensen's coefficient $M_K(h)$ exists.

\begin{theorem} \label{thm1}
Let $A$ be a self-adjoint operator on
a Hilbert space $\mathcal{H}$ and assume that $Sp(A)\subseteq [m,M]$ for some scalars $m,M$ with
$0<m < M$.
Let $h:[0,1]\to\Bbb{R}$ be a non-negative and non-zero function.
If $f$ is a continuous $h$-convex function on $[m, M]$, then
\begin{equation}\label{eq0-thm1}
f(\langle Ax, x\rangle) \leq M_{(0,1)}(h)\langle f(A)x,x\rangle
\end{equation}
for each $x\in \mathcal{H}$ with $||x||=1$.
\end{theorem}
\begin{proof}
It is clear that if $M_{(0,1)}(h)=+\infty$, then the inequality \eqref{eq0-thm1} holds.
Assume that $M_{(0,1)}(h)<+\infty$.
It follows from $h$-convexity of $f$ that
\begin{equation}\label{eq1-thm1}
f(\lambda a+(1-\lambda)b)-h(1-\lambda)f(b)\leq h(\lambda) f(a)
\end{equation}
for all $a,b\in[m,M]$ and for all $\lambda\in (0,1)$.
By dividing both sides of \eqref{eq1-thm1} with $\lambda\in (0,1)$, one can reach
\begin{equation}\label{eq1-eq1-thm1}
\frac{f(\lambda a+(1-\lambda)b)-h(1-\lambda)f(b)}{\lambda}\leq \frac{h(\lambda)}{\lambda} f(a)
\end{equation}
for all $a,b\in[m,M]$.
Define
\begin{equation}\label{eq2-thm1}
\alpha:=\min_{b\in[m,M]}\frac{f(\lambda a+(1-\lambda)b)-h(1-\lambda)f(b)}{\lambda(a-b)}.
\end{equation}
The inequalities \eqref{eq1-eq1-thm1} and \eqref{eq2-thm1} entail that
$$
\alpha(a-b)\leq \frac{h(\lambda)}{\lambda}f(a)
$$
for all $a,b\in[m,M]$ and $\lambda\in (0,1)$.
Consider the linear function $l(t):=\alpha(t-b)$ and $\bar{g}=\langle Ax,x\rangle$.
This implies that $l(a)\leq \frac{h(\lambda)}{\lambda}f(a)$ for all $a\in[m,M]$ and
$m\leq \bar{g}\leq M$, respectively.
Consider the straight line $l'(t):=\alpha(t-\bar{g})+f(\bar{g})$ passing through the point
$(\bar{g},f(\bar{g}))$ and parallel to the line $l$.
The continuity of the function $f$ ensures that
\begin{equation}\label{eq3-thm1}
l'(\bar{g})\geq f(\bar{g})-\epsilon
\end{equation}
for all $\epsilon > 0$.
We now consider two cases:

(i) Assume that $l'(t)\leq \frac{h(\lambda)}{\lambda}f(t)$ for every $t\in[m,M]$.
By using the functional calculus, one has
$l'(A)\leq \frac{h(\lambda)}{\lambda}f(A)$ and consequently
\begin{equation}\label{eq4-thm1}
\langle l'(A)x,x \rangle\leq \frac{h(\lambda)}{\lambda}\langle f(A)x,x \rangle
\end{equation}
for each $x\in \mathcal{H}$ with $||x||=1$.
The linearity of the function $l'$ and the inequalities \eqref{eq3-thm1} and \eqref{eq4-thm1} imply
\begin{align*}%\label{eq5-thm1}
f(\langle Ax,x\rangle)-\epsilon\leq l'(\langle Ax,x\rangle)=\langle l'(A)x,x\rangle\leq \frac{h(\lambda)}{\lambda}\langle f(A)x,x \rangle.
\end{align*}
Since $\epsilon$ is arbitrary, we observe that
\begin{equation}\label{eq5-thm1}
f(\langle Ax,x\rangle)\leq \frac{h(\lambda)}{\lambda}\langle f(A)x,x \rangle.
\end{equation}

(ii) Assume that there exist some points $t\in[m,M]$ such that $l'(t)>\frac{h(\lambda)}{\lambda}f(t)$.
Define the sets $T$ and $S$ as follows:
$$
T:=\{t\in [m,\bar{g}]: l'(t)>\frac{h(\lambda)}{\lambda}f(t)\},
$$
$$
S:=\{t\in [\bar{g},M]: l'(t)>\frac{h(\lambda)}{\lambda}f(t)\}.
$$
Consider $t_T:=\max\{t: t\in T\}$ and $t_S:=\min\{t: t\in S\}$.
We use two lines passing through the points $(t_T,0)$, $(\bar{g},f(\bar{g}))$
and $(t_S,0)$ and $(\bar{g},f(\bar{g}))$, respectively.
Let $l_T$ be the line passing through the points $(t_T,0)$ and $(\bar{g},f(\bar{g}))$
and $l_S$ the line passing through the points $(t_S,0)$ and $(\bar{g},f(\bar{g}))$
and define the function $L$ as follows:
$$
L(t):=\Big\{
\begin{array}{ll}
l_T(t), t\in [m,\bar{g}],\\
l_S(t), t\in [\bar{g},M].
\end{array}
$$
We prove that the inequality $L(t)\leq \frac{h(\lambda)}{\lambda}f(t)$ holds for all $t\in[m, M]$.
We consider the partition $\{m, t_T, \bar{g}, t_S, M\}$ for the closed interval $[m, M]$
and remark that $l_T(t)\leq 0$ for every $t\in[m, t_T]$. Since $f(t)\geq 0$, we clearly observe that
$l_T(t)\leq \frac{h(\lambda)}{\lambda}f(t)$ for every $t\in[m, t_T]$.
On the other hand, we see that
\begin{equation}\label{main1}
l'(t)\leq \frac{h(\lambda)}{\lambda}f(t)
\end{equation}
for every $t\in (t_T, \bar{g}]$; otherwise,
there exists $t_0\in (t_T, \bar{g}]$ such that $l'(t_0)> \frac{h(\lambda)}{\lambda}f(t_0)$ and
so $t_0\in T$ and $t_0< t_T$, which is a contradiction.
So, it follows from \eqref{main1} that
\begin{equation}\label{main110}
l'(t_T)\leq \frac{h(\lambda)}{\lambda}f(t_T)
\end{equation}
by letting $t$ tends to $t_T$ from right in \eqref{main1}.
Moreover, since $t_T\in \bar{T}$, the reverse inequality holds in \eqref{main110} and hence $l'(t_T)=\frac{h(\lambda)}{\lambda}f(t_T)$.
It follows that $l'$ is the line passing through the points $(t_T,\frac{h(\lambda)}{\lambda}f(t_T))$ and $(\bar{g},f(\bar{g}))$
and its slope is $\alpha=\frac{f(\bar{g})-\frac{h(\lambda)}{\lambda}f(t_T)}{\bar{g}-t_T}$, where
the slope of $l_T$ is $\alpha'=\frac{f(\bar{g})}{\bar{g}-t_T}$.
By the inequality \eqref{main1} we observe that
$$
l_T(t)=\alpha'(t-\bar{g})+f(\bar{g})\leq \alpha(t-\bar{g})+f(\bar{g})=l'(t)\leq \frac{h(\lambda)}{\lambda}f(t)
$$
for every $t\in (t_T, \bar{g}]$.
So, $L(t)=l_T(t)\leq\frac{h(\lambda)}{\lambda}f(t)$ for every $t\in [m, \bar{g}]$.

By the similar methods one can show that $L(t)=l_S(t)\leq\frac{h(\lambda)}{\lambda}f(t)$ for every $t\in [\bar{g}, M]$.
Note that
the lines $l_T$ and $l_S$ are joining at the point along the length of $\bar{g}$ and so
$l_T(\bar{g})=l_S(\bar{g})$ and since $f$ is continuous,
\begin{equation}%\label{eq5-thm1}
l_T(\bar{g})=f(\bar{g})\geq f(\bar{g})-\epsilon
\end{equation}
for arbitrary $\epsilon > 0$. For the case  $Sp(A)\subseteq [m,\bar{g}]$, we have
\begin{align*}%\label{eq5-thm1}
f(\langle Ax,x\rangle)-\epsilon\leq l_T(\langle Ax,x\rangle)=\langle l_T(A)x,x\rangle\leq \frac{h(\lambda)}{\lambda}\langle f(A)x,x \rangle.
\end{align*}
Moreover, for the case $Sp(A)\subseteq [\bar{g},M]$, we have
\begin{align*}%\label{eq5-thm1}
f(\langle Ax,x\rangle)-\epsilon
&\leq l_T(\langle Ax,x\rangle)=l_S(\langle Ax,x\rangle)=\langle l_S(A)x,x\rangle\\
&\leq \frac{h(\lambda)}{\lambda}\langle f(A)x,x \rangle
\end{align*}
and consequently we can deduce \eqref{eq5-thm1}.
Now, by taking the infimum over all $\lambda\in (0,1)$ on both sides of the inequality \eqref{eq5-thm1},
one can deduce \eqref{eq0-thm1}.
\end{proof}

%^^^^^^^^^^^^^^^^^^^^
\begin{corollary}\label{thebest}
Let the conditions of Theorem \ref{thm1} be satisfied.
If the function $\frac{h(t)}{t}$ is decreasing on $(0,1)$, then $\lambda=\frac{1}{2}$
is the best possible in the inequality \eqref{eq5-thm1}. Indeed, one has
\begin{equation}\label{rem-eq}
f(\langle Ax, x\rangle) \leq 2h(\frac{1}{2})\langle f(A)x,x\rangle
\end{equation}
for each $x\in \mathcal{H}$ with $||x||=1$.
\end{corollary}
\begin{proof}
We claim that $\frac{1}{2}$ is the best possible for $\lambda$ in the inequality \eqref{eq5-thm1}.
We divide the interval $(0,1)$ to two parts $(0,\frac{1}{2}]$ and $(\frac{1}{2},1)$, respectively.
In the first part by using the inequality \eqref{eq5-thm1} we detect the best possible % for $\lambda\in (0,\frac{1}{2}]$.
and in the second part by an example we demonstrate that the inequality \eqref{eq5-thm1} does not hold
in general. % for $\lambda\in (\frac{1}{2},1)$.
In fact, in this situation we conclude that the inequality \eqref{eq5-thm1} only holds in the first part. %for all $\lambda\in (0,\frac{1}{2}]$.

(1) Let $0<\lambda\leq\frac{1}{2}$. Since the function $\frac{h(t)}{t}$ is decreasing,
 its infimum value over $(0,\frac{1}{2}]$ occurs at the endpoint $\lambda=\frac{1}{2}$
 and so $M_{(0,\frac{1}{2}]}(h)=\frac{h(\frac{1}{2})}{\frac{1}{2}}=2h(\frac{1}{2})$.

(2) Let $\frac{1}{2}<\lambda<1$.
We show that there is an $h$-convex function such that does not satisfy the inequality \eqref{eq5-thm1}.
Let $h(t)=\sqrt{t}$ %, $0<\beta<\frac{1}{2}$,
and $t>0$.
Then, $\frac{h(t)}{t}=t^{-\frac{1}{2}}$ is decreasing on $(0,1)$.
Define $g: [0, \infty)\to \Bbb{R}$ by $g(t)=\sqrt{t}$.
%$$
%g(t)=\bigg\{
%\begin{array}{ll}
%1, & t\neq\frac{a}{2}\\
%2^{1-k}, & t=\frac{a}{2}.
%\end{array}
%$$
Note that $g$ is $h$-convex, since % \cite[Example 6]{S.Var-2007}.
\begin{align*}
g(\alpha x+(1-\alpha)y)&=(\alpha x+(1-\alpha)y)^{\frac{1}{2}}\\
&\leq (\alpha x)^{\frac{1}{2}}+((1-\alpha)y)^{\frac{1}{2}}\\
&=\alpha^{\frac{1}{2}} x^{\frac{1}{2}}+(1-\alpha)^{\frac{1}{2}}y^{\frac{1}{2}}\\
%&\leq \alpha^{\beta} x^{\frac{1}{2}}+(1-\alpha)^{\beta}y^{\frac{1}{2}}\\
&=h(\alpha)x^{\frac{1}{2}}+h(1-\alpha)y^{\frac{1}{2}}
\end{align*}
for every $x,y\geq 0$ and $\alpha\in [0,1]$.
Consider
$
A=\Big(
\begin{array}{ll}
1 & 0\\
0 & 0
\end{array}
 \Big)
$
and $x=(\frac{1}{\sqrt{2}},\frac{1}{\sqrt{2}})$.
A simple calculation shows that $g(\langle Ax, x\rangle)=g(\frac{1}{2})=\sqrt{\frac{1}{2}}$ and $\langle g(A)x,x\rangle=\frac{1}{2}$.
%We have $\lambda^{\beta-1}\leq \lambda^{-1}$, since $-1<\beta-1$.
%Note that $\frac{1}{2}<\lambda$ implies $\lambda^{-1}<2$.
%This indicates that $\lambda^{\beta-1}<2$.
Since $g$ is $h$-convex, by \eqref{eq5-thm1}, we have
$$
g(\langle Ax, x\rangle) \leq \frac{h(\lambda)}{\lambda}\langle g(A)x,x\rangle
$$
and so $\sqrt{\frac{1}{2}}\leq\frac{\lambda^{-\frac{1}{2}}}{2}$
which is a contradiction.
\end{proof}

%such that $h(\lambda)\geq \lambda$ for every $\lambda\in [0,1]$
%%%%%%%%%%%%%%%%%%%%%%%%%%%%%%%%%%%%
%%%%%%%%%%%%%%%%%%%%%%%%%%%%%%%%%%%%

%\begin{corollary} \label{cor1}
%Under the hypotheses of Theorem \ref{thm1},
%if $x \in \mathcal{H}$, $||x||\neq 1$, then %$\langle x,x\rangle \neq 1$
%\begin{equation}\label{eq0-cor1}
%f\Big(\frac{\langle Ax, x\rangle}{\langle x, x\rangle}\Big) \leq 2h\bigg(\frac{1}{2}\bigg)\frac{\langle f(A)x,x\rangle}{\langle x, x\rangle}.
%\end{equation}
%\end{corollary}
%\begin{proof}
%Let $y:=\frac{x}{\sqrt{\langle x, x\rangle}}$ and apply Theorem \ref{thm1}.
%\end{proof}
%%%%%%%%%%%%%%%%%%%%%%%%%%%%%%%
%%%%%%%%%%%%%%%%%%%%%%%%%%%%%%%

\begin{corollary} \label{cor-part1}
Let $A$ be a self-adjoint operator on
a Hilbert space $\mathcal{H}$ and assume that $Sp(A)\subseteq [m,M]$ for some scalars $m,M$ with
$0<m < M$ and $x\in \mathcal{H}$ with $||x||=1$.
\begin{itemize}
\item [(1)] If $f$ is a continuous convex function on $[m, M]$, then
\begin{equation}
f(\langle Ax, x\rangle)\leq \langle f(A)x,x\rangle.
\end{equation}
\item [(2)] If $f$ is a continuous $P$-class function on $[m, M]$, then
\begin{equation}
f(\langle Ax, x\rangle)\leq 2\langle f(A)x,x\rangle.
\end{equation}
\item [(3)] If $f$ is a continuous $s$-convex function on $[m, M]$ in the second sense, then
\begin{equation}
f(\langle Ax, x\rangle)\leq 2^{1-s}\langle f(A)x,x\rangle.
\end{equation}
\item [(4)] If $f$ is a continuous Godunova-Levin function on $[m, M]$, then
\begin{equation}
f(\langle Ax, x\rangle)\leq 4\langle f(A)x,x\rangle.
\end{equation}
\item [(5)] If $f$ is a continuous $s$-Godunova-Levin function on $[m, M]$, then
\begin{equation}
f(\langle Ax, x\rangle)\leq 2^{1+s}\langle f(A)x,x\rangle.
\end{equation}
\end{itemize}
\end{corollary}
\begin{proof}
Consider $h(t)=t$, $h(t)=1$, $h(t)=t^s$, $h(t)=\frac{1}{t}$, and $h(t)=\frac{1}{t^s}$ in parts (1)-(5), respectively.
Then, the function $\frac{h(t)}{t}$ is decreasing in each part and Corollary \ref{thebest}
implies the Jensen's coefficient is $2h(\frac{1}{2})$.
So, a simple calculation gets the desired result in each part.
\end{proof}

%-------------------------------------------------------------------------------------------------------
%
%
%%%%%%%%%%%%%%%%%%%%%%%%%%%%%%%%%%%%%%%%%%%%%%%%%%%%%%%%%%%%%%%%%%%%%%%%%%%%%%%%%%%%%%%%%%%%%%%%%%%%%%%%

We provide a refinement of the Mond and Pe$\check{c}$ari$\acute{c}$ inequality for $h$-convex functions.

\begin{corollary}\label{Es-cor-refine}
Let the conditions of Theorem \ref{thm1} be satisfied.
If
%$f$ is a continuous convex function on $[m, M]$ such that
%$f(\langle Ax, x\rangle) < \langle f(A)x,x\rangle$ and
$M_{(0,1)}(h)<1$ and $f(\langle Ax, x\rangle) < \langle f(A)x,x\rangle$, then
\begin{equation}%\label{eq0-thm1}
f(\langle Ax, x\rangle) \leq M_{(0,1)}(h)\langle f(A)x,x\rangle < \langle f(A)x,x\rangle
\end{equation}
for each $x\in \mathcal{H}$ with $||x||=1$. %and $A\in \Gamma(f,m,x)$.
\end{corollary}
%%%%%%%%%%%%%%%%%%%%%%%%%%%%%%%
\begin{proof}
The first inequality follows from Theorem \ref{thm1} and the second one follows
from the fact that $M_{(0,1)}(h)<1$.
\end{proof}
We remark that in Corollary \ref{Es-cor-refine} the condition
$f(\langle Ax, x\rangle) < \langle f(A)x,x\rangle$
is essential, since we reach the contradiction $M_{(0,1)}(h)=1$ when equality holds.

%%%%%%%%%%%%%%%%%%%%%%%%%%%%%%%
%%%%%%%%%%%%%%%%%%%%%%%%%%%%%%%
%%%%%%%%%%%%%%%%%%%%%%%%%%%%%%%
%%%%%%%%%%%%%%%%%%%%%%%%%%%%%%%
%%%%%%%%%%%%%%%%%%%%%%%%%%%%%%%
%
\begin{theorem}\label{thm3}
Let the conditions of Theorem \ref{thm1} be satisfied. Then,
\begin{equation}\label{eq0-thm3}
\langle f(A)x, x\rangle\leq M_{(0,1)}(h)\Big(\frac{M-\langle Ax,x\rangle}{M-m}f(m)+\frac{\langle Ax,x\rangle-m}{M-m}f(M)\Big).
\end{equation}
\end{theorem}
%%%%%%%%%%%%%%%%%%%%%%%%%%%%%%%
\begin{proof}
Consider
$D=\Big(
\begin{array}{ll}
m & 0\\
0 & M
\end{array}
 \Big)$
 and
 $x=\Bigg(
\begin{array}{l}
\sqrt{\frac{M-t}{M-m}}\\
\sqrt{\frac{t-m}{M-m}}
\end{array}
 \Bigg)$.
By applying \ref{thm1} we have
\begin{align*}
f(t)&=f(\langle Dx,x \rangle)\\
&\leq M_{(0,1)}(h)\langle f(D)x,x\rangle\\
&=M_{(0,1)}(h)\Big(\frac{M-t}{M-m}f(m)+\frac{t-m}{M-m}f(M)\Big).
\end{align*}
Since the operator $M_{(0,1)}(h)\Big(\frac{M-A}{M-m}f(m)+\frac{A-m}{M-m}f(M)\Big)-f(A)$ is positive, we get \eqref{eq0-thm3}.
\end{proof}

%%%%%%%%%%%%%%%%%%%%%%%%%%%%%%%%%%%%%%%%%%%%%%%%%%%%%%%%%
%
\begin{theorem}\label{thm4}
Let the conditions of Theorem \ref{thm1} be satisfied.
Let $J$ be an interval such that $f([m, M])\subset J$.
If $F(u, v)$ is a real function defined on $J\times J$ and non--decreasing in $u$, then
\begin{align}\label{eq0-thm4}
F&(\langle f(A)x,x\rangle,f(\langle Ax,x\rangle))\nonumber\\
&\leq\max_{t\in[m,M]}F\bigg(M_{(0,1)}(h)\bigg(\frac{M-t}{M-m}f(m)+\frac{t-m}{M-m}f(M)\bigg),f(t)\bigg)\nonumber\\
&=\max_{\theta\in[0,1]}F\bigg(M_{(0,1)}(h)(\theta f(m)+(1-\theta)f(M)),f(\theta m+(1-\theta)M)\bigg).
\end{align}
\end{theorem}
%%%%%%%%%%%%%%%%%%%%%%%%%%%%%%%
\begin{proof}
Since $\bar{g} = \langle Ax,x\rangle\in [m, M]$, by the non-decreasing character of $F$ and Theorem \ref{thm3}, one has
\begin{align*}
F&(\langle f(A)x,x\rangle,f(\langle Ax,x\rangle))\nonumber\\
&\leq F\bigg(M_{(0,1)}(h)\Big(\frac{M-\bar{g}}{M-m}f(m)+\frac{\bar{g}-m}{M-m}f(M)\Big),f(\bar{g})\bigg)\nonumber\\
&\leq \max_{t\in[m,M]}F\bigg(M_{(0,1)}(h)\bigg(\frac{M-t}{M-m}f(m)+\frac{t-m}{M-m}f(M)\bigg),f(t)\bigg).
\end{align*}
The second form of the right side of \eqref{eq0-thm4} follows
at once from the change of variable $\theta=\frac{M-t}{M-m}$, so that
$t=\theta m+(1-\theta)M$, with $0\leq \theta\leq 1$.
\end{proof}

%More simply just by replacing $F$ by $-F$ in the above theorem one can prove the reverse inequality as follows.
%%%%%%%%%%%%%%%%%%%%%%%%%%%%%%%%%%%%%%%%%%%%%%%%%%%%%%%%%
%
%\begin{corollary}\label{cor3}
%Under the same hypotheses as Theorem \ref{thm4}, except that $F$ is
%non--increasing in its first variable, we have
%\begin{align*}%\label{eq0-cor3}
%F&(\langle f(A)x,x\rangle,f(\langle Ax,x\rangle))\nonumber\\
%&\geq\min_{t\in[m,M]}F\bigg(M_{(0,1)}(h)\bigg(\frac{M-t}{M-m}f(m)+\frac{t-m}{M-m}f(M)\bigg),f(t)\bigg)\nonumber\\
%&=\min_{\theta\in[0,1]}F\bigg(M_{(0,1)}(h)(\theta f(m)+(1-\theta)f(M)),f(\theta m+(1-\theta)M)\bigg).
%\end{align*}
%\end{corollary}

%%%%%%%%%%%%%%%%%%%%%%%%%%%%%%%%%%%%%%%%%%%%%%%%%%%%%%%%%%%%%%%%%%%%%%%%%%%%%%%%%%%%%%%%%%%%%%%

\begin{definition}\label{PCTD}
The function $f$ is piecewise continuously twice differentiable on $[m,M]$ whenever the following conditions fulfil:
\begin{itemize}
\item[(1)] $f$ is continuous on $[m,M]$,
\item[(2)] there exists a finite subdivision $\{x_0,...,x_n\}$ of $[m,M]$, $x_0=a$, $x_n=b$ such that
\subitem (2.1) $f$ is continuously twice differentiable on $(x_{i-1},x_i)$ for every $i\in \{1, ...,n\}$,
\subitem (2.2) the one-sided limits $\lim_{x\to x_{i-1}^+}f'(x)$ and $\lim_{x\to x_i^-}f'(x)$ exist for every $i\in\{1, ...,n\}$,
\subitem (2.3) the one-sided limits $\lim_{x\to x_{i-1}^+}f''(x)$ and $\lim_{x\to x_i^-}f''(x)$ exist for every $i\in\{1, ...,n\}$.
\end{itemize}
\end{definition}

%Note that the function $f$ defined in Example \ref{exam} is piecewise continuously twice differentiable on $[0,\infty)$.

We provide a converse inequality in Theorem \ref{thm1}.
\begin{theorem} \label{cnrs-thm4}
Let the conditions of Theorem \ref{thm1} be satisfied.
Moreover, let $f$ be piecewise continuously twice differentiable
on $[m,M]$.
\begin{itemize}
\item[(i)] There exists $\alpha>0$ such that
$$
\frac{1}{M_{(0,1)}(h)\alpha}\langle f(A)x,x\rangle\leq f(\langle Ax, x\rangle).
$$
\item[(ii)] There exists $\beta>0$ such that
$$
\frac{1}{M_{(0,1)}(h)}\langle f(A)x,x\rangle-\beta\leq f(\langle Ax, x\rangle).
$$
\end{itemize}
\end{theorem}
%%%%%%%%%%%%%%%%%%%%%%%%%%%%%%%
\begin{proof}
(i) Suppose $R=\{x_0, x_1, ..., x_{n}\}$, $x_0=m$, $x_n=M$ is the finite subdivision of $[m,M]$ such that
the conditions of Definition \ref{PCTD} fulfils.
Consider $F(u,v)=\frac{u}{v}$, $J=(0,\infty)$, $\varphi_h(t)=M_{(0,1)}(h)\varphi_i(t)$ for every $t\in [x_{i-1}, x_i]$,
where $\varphi_i(t)=\frac{L_i(t)}{f(t)}$, $L_i(t)=f(x_{i-1})+\mu_i(t-x_{i-1})$ and $\mu_i=\frac{f(x_i)-f(x_{i-1})}{x_i-x_{i-1}}$.
According Theorem \ref{thm4} we have
\begin{align}\label{max-m}
\frac{\langle f(A)x,x\rangle}{f(\langle Ax,x\rangle)}\leq \max_{t\in[m,M]}\varphi_h(t)= M_{(0,1)}(h)\max_{1\leq i\leq n}\max_{t\in[x_{i-1}, x_i]}\varphi_i(t).
\end{align}
Now $\varphi_i'(t)= \frac{G_i(t)}{f(t)^2}$, where $G_i(t)=\mu_i f(t) - L_i(t)f'(t)$ for every $t\in [x_{i-1}, x_i]$.
If $\mu_i\neq 0$, then
$$
G_i'(t)=-L_i(t)f''(t).
$$
If $\mu_i = 0$ and $\bar{t_i}\in (x_{i-1}, x_i)$ is the unique solution of the equation $f'(t) = 0$, then
we consider
$$
\lambda_i = \max_{t\in [x_{i-1},x_{i}]} \varphi_h(t)=M_{(0,1)}(h)\frac{f(x_{i-1})}{f(\bar{t_i})}.
$$
Define
$$
A=\{i: \lim_{t\to x_{i-1}^+}f''(t)>0, \lim_{t\to x_{i}^-}f''(t)>0, f''(t)>0, t\in (x_{i-1}, x_i)\},
$$
$$
B=\{i: \lim_{t\to x_{i-1}^+}f''(t)<0, \lim_{t\to x_{i}^-}f''(t)<0, f''(t)<0, t\in (x_{i-1}, x_i)\}.
$$
Suppose $t\in [m,M]$. Then, there exists $i\in \{1, ..., n\}$ such that $t\in [x_{i-1}, x_i]$.

(1) If $i\in A$, then $G_i'(t)<0$ and so $G_i$ is decreasing on $[x_{i-1}, x_i]$.
So,
$$
G_i(x_{i-1})G_i(x_i)=-f(x_{i-1})f(x_i)(\mu_i-f'(x_{i-1}))(f'(x_i)-\mu_i)<0,
$$
This indicates  the equation $G_i(t)=0$ has a unique solution at $\bar{t_i}\in (x_{i-1}, x_i)$
and so the equation $\varphi_i'(t)=0$ has a unique solution at $\bar{t_i}\in (x_{i-1}, x_i)$.
Let
$D_i=\Big(
\begin{array}{ll}
x_{i-1} & 0\\
0 & x_i
\end{array}
 \Big)$
 and
 $x=\Bigg(
\begin{array}{l}
\sqrt{\frac{x_i-t}{x_i-x_{i-1}}}\\
\sqrt{\frac{t-x_{i-1}}{x_{i}-x_{i-1}}}
\end{array}
 \Bigg).$
Since $i\in A$, the function $f$ is convex on $[x_{i-1},x_{i}]$.
So, by Theorem \ref{thm0}, for the convex function $f$ on $[x_{i-1},x_{i}]$, one can reach
\begin{align*}
f(t)&=f(\langle D_ix,x \rangle)\\
&\leq \langle f(D_i)x,x\rangle=\frac{x_{i}-t}{x_{i}-x_{i-1}}f(x_{i-1})+\frac{t-x_{i-1}}{x_{i}-x_{i-1}}f(x_{i})=L_i(t).
\end{align*}
Consequently, $\frac{L_i(t)}{f(t)}\geq 1$ for every $t\in [x_{i-1}, x_i]$ and
$$
\varphi_h(t)=M_{(0,1)}(h)\frac{L_i(t)}{f(t)}\geq M_{(0,1)}(h)
$$
for every $t\in [x_{i-1}, x_i]$ where the equality occurs at $x_{i-1}$ and $x_i$.
Note that the maximum value of $\varphi_i$ is attained in $\bar{t_i}\in [x_{i-1}, x_i]$,
since $\varphi_i''(\bar{t_i})<0$.
We consider
\begin{align*}
\lambda_i&=\max_{t\in [x_{i-1}, x_i]}\varphi_h(t)\\
&=M_{(0,1)}(h)\max_{t\in [x_{i-1}, x_i]}\varphi_i(t)\\
&=M_{(0,1)}(h)\varphi_i(\bar{t_i})\\
&=M_{(0,1)}(h)\frac{L_i(\bar{t_i})}{f(\bar{t_i})}\\
&=M_{(0,1)}(h)\frac{\mu_i}{f'(\bar{t_i})}.
\end{align*}
The last equality comes from the fact that $G_i(\bar{t_i})=0$.

(2) If $i\in B$, then define
$D_i$ and  $x$ as the part (1) and apply Theorem \ref{thm0} for the concave function $f$ on $[x_{i-1},x_{i}]$.
So, $\frac{L_i(t)}{f(t)}\leq 1$ and this inequality yields
$$
0\leq \varphi_h(t)=M_{(0,1)}(h)\frac{L_i(t)}{f(t)}\leq M_{(0,1)}(h)
$$
for every $t\in [x_{i-1}, x_i]$ where equality occurs at $x_{i-1}$ and $x_i$. We consider
$$
\lambda_i = \max_{t\in [x_{i-1},x_{i}]} \varphi_h(t) = M_{(0,1)}(h).
$$
It follows from the cases (1) and (2) in the part (i) that
$$
\lambda_i =
\Bigg\{
\begin{array}{lll}
M_{(0,1)}(h)\frac{f(x_{i-1})}{f(\bar{t_i})}, & \mu_i=0, i\in \{1, ..., n\}\backslash A\cup B,\\
M_{(0,1)}(h)\frac{\mu_i}{f'(\bar{t_i})}, & \mu_i\neq 0, i\in A,\\
M_{(0,1)}(h), & \mu_i\neq 0, i\in B.
\end{array}
$$
Define $\lambda=\max_{1\leq i\leq n}\lambda_i$.
Then, $\lambda=M_{(0,1)}(h)\alpha$, where
$$
\alpha=\max\{\max_{i\in(A\cup B)^c}\frac{f(x_{i-1})}{f(\bar{t_i})}, \max_{i\in A}\frac{\mu_i}{f'(\bar{t_i})}, 1\}.
$$
By virtue of \eqref{max-m}, we deduce
$$
\frac{\langle f(A)x,x\rangle}{f(\langle Ax,x\rangle)}\leq \max_{t\in[m,M]}\varphi_h(t)= M_{(0,1)}(h)\alpha.
$$
%If the sign of $f''$ is negative, then $A=\emptyset$ and so $\lambda=\max_{i\in B}\lambda_i=2^{1-s}$. Hence,
%$$
%\frac{\langle f(A)x,x\rangle}{f(\langle Ax,x\rangle)}\leq \max_{t\in[m,M]}\varphi_h(t)= \lambda=2^{1-s}.
%$$

(ii)
Consider the sets $R$, $P$, $A$, and $B$ as the part (i)
and define $F(u,v)=u-M_{(0,1)}(h)v$, $J=\Bbb{R}$ and $\varphi_h(t)=M_{(0,1)}(h)\varphi_i(t)$ for every $t\in [x_{i-1}, x_i]$,
where $\varphi_i(t)=L_i(t)-f(t)$ and $L_i(t)$ defined in the part (i).
By virtue of Theorem \ref{thm4} we yield
\begin{align}\label{max-m2}
\langle f(A)x,x\rangle - M_{(0,1)}(h)f(\langle Ax,x\rangle)&\leq \max_{t\in[m,M]}\varphi_h(t)\nonumber\\
&= M_{(0,1)}(h)\max_{1\leq i\leq n}\max_{t\in[x_{i-1}, x_i]}\varphi_i(t).
\end{align}
If $\mu_i\neq 0$, then
$\varphi_i''(t)=-f''(t)$ and if $\mu_i = 0$ and $\bar{t_i}\in (x_{i-1},x_i)$ is the unique solution of the equation $f'(t) = 0$, then
we define
$$
\lambda_i = \max_{t\in [x_{i-1}, x_i]}\varphi_h(t)=M_{(0,1)}(h)(f(x_{i-1})-f(\bar{t_i})).
$$
Suppose $t\in [m,M]$. Then, there exists $i\in \{1, ..., n\}$ such that $t\in [x_{i-1}, x_i]$.

(1) If $i\in A$, then $\varphi_i''(t)<0$ for every $t\in[x_{i-1}, x_i]$ and so $\varphi_i'$ is decreasing on $[x_{i-1}, x_i]$.
On the other hand,
the equation $\varphi_i'(t)=0$ has a unique solution at $t=\bar{t_i}\in [x_{i-1}, x_i]$, since $\varphi_i'(x_{i-1})\varphi_i'(x_i)<0$.
Clearly, $\varphi_i''(\bar{t_i})<0$ and so %$\bar{t_i}$ is the global  $\varphi_i$
the maximum value of $\varphi_i$ is attained in $\bar{t_i}$.
We define
\begin{align*}
\lambda_i&=\max_{t\in [x_{i-1}, x_i]}\varphi_h(t)\\
&=M_{(0,1)}(h)\max_{t\in [x_{i-1}, x_i]}\varphi_i(t)\\
&=M_{(0,1)}(h)\varphi_i(\bar{t_i})\\
&=M_{(0,1)}(h)(L_i(\bar{t_i})-f(\bar{t_i}))\\
&=M_{(0,1)}(h)(f(x_{i-1})+\mu_i(\bar{t_i}-x_{i-1})-f(\bar{t_i})).
\end{align*}

(2) If $i\in B$, then $f''(t)<0$ for every $t\in[x_{i-1}, x_i]$. This means that $f$ is concave on $[x_{i-1}, x_i]$
and so $f(t)\geq L_i(t)$ for every $t\in[x_{i-1}, x_i]$.
This ensures $\varphi_i(t)=L_i(t)-f(t)\leq 0$ and this inequality entails
$$
\max_{t\in[x_{i-1}, x_i]}\varphi_i(t)\leq 0.
$$
Since $\varphi_i(x_i)=0=\varphi_i(x_{i-1})$, $\varphi_i$ attains its maximum value and the maximum value is $0$.
So that
$$
\lambda_i = \max_{t\in [x_{i-1},x_{i}]} \varphi_h(t) = M_{(0,1)}(h)\max_{t\in[x_{i-1}, x_i]}\varphi_i(t)=\varphi_i(x_i)=0.
$$
Consequently, it follows from the cases (1) and (2) in the part (ii) that
$$
\lambda_i = \Bigg\{ \begin{array}{lll}
M_{(0,1)}(h)(f(x_{i-1})-f(\bar{t_i})), & \mu_i=0, i\in \{1, ..., n\}\backslash A\cup B,\\
M_{(0,1)}(h)(f(x_{i-1})+\mu_i(\bar{t_i}-x_{i-1})-f(\bar{t_i})), & \mu_i\neq 0, i\in A,\\
0, &  \mu_i\neq 0, i\in B.
\end{array}
$$
Define $\lambda=\max_{1\leq i\leq n}\lambda_i$.
Then, $\lambda=M_{(0,1)}(h)\beta$, where
$$
\beta=\max\{\max_{i\in(A\cup B)^c}(f(x_{i-1})-f(\bar{t_i})), \max_{i\in A}(f(x_{i-1})+\mu_i(\bar{t_i}-x_{i-1})-f(\bar{t_i})), 0\}.
$$
In view of \eqref{max-m2}, we deduce
$$
\langle f(A)x,x\rangle-M_{(0,1)}(h)f(\langle Ax,x\rangle)\leq \max_{t\in[m,M]}\varphi_h(t)=\lambda=M_{(0,1)}(h)\beta.
$$
\end{proof}

% %%%%%%%%%%%%%%%%%%%%%%%%%%%%%%%%%%%%%%%%%%%%%%%%%%%%%%%%
%%%%%%%%%%%%%%%%%%%%%%%%%%%%%%%
%%%%%%%%%%%%%%%%%%%%%%%%%%%%%%%
%\begin{itemize}
%\item[(i)] There exists $\alpha>0$ such that
%$$
%\frac{1}{M_{(0,1)}(h)\alpha}\langle f(A)x,x\rangle\leq f(\langle Ax, x\rangle).
%$$
%\item[(ii)] There exists $\beta>0$ such that
%$$
%\frac{1}{M_{(0,1)}(h)}\langle f(A)x,x\rangle-\beta\leq f(\langle Ax, x\rangle).
%$$
%\end{itemize}

% %%%%%%%%%%%%%%%%%%%%%%%%%%%%%%%%%%%%%%%%%%%%%%%%%%%%%%%%
%%%%%%%%%%%%%%%%%%%%%%%%%%%%%%%
%%%%%%%%%%%%%%%%%%%%%%%%%%%%%%%

\begin{corollary} \label{cor-part2}
Let the function $f$ be a piecewise continuously twice differentiable on $[m, M]$ and $A$ a self-adjoint operator on
a Hilbert space $\mathcal{H}$.
Assume that $Sp(A)\subseteq [m,M]$ for some scalars $m,M$ with
$0<m < M$ and $x\in \mathcal{H}$ with $||x||=1$.
\begin{itemize}
\item [(1)] If $f$ is convex on $[m, M]$, then
\begin{itemize}
\item[(i)] there exists $\alpha>0$ such that
$$
\frac{1}{\alpha}\langle f(A)x,x\rangle\leq f(\langle Ax, x\rangle),
$$
\item[(ii)] there exists $\beta>0$ such that
$$
\langle f(A)x,x\rangle-\beta\leq f(\langle Ax, x\rangle).
$$
\end{itemize}
\item [(2)] If $f$ is $P$-class on $[m, M]$, then
\begin{itemize}
\item[(i)] there exists $\alpha>0$ such that
$$
\frac{1}{2\alpha}\langle f(A)x,x\rangle\leq f(\langle Ax, x\rangle),
$$
\item[(ii)] there exists $\beta>0$ such that
$$
\frac{1}{2}\langle f(A)x,x\rangle-\beta\leq f(\langle Ax, x\rangle).
$$
\end{itemize}
\item [(3)] If $f$ is $s$-convex on $[m, M]$ in the second sense, then
\begin{itemize}
\item[(i)] there exists $\alpha>0$ such that
$$
\frac{1}{2^{1-s}\alpha}\langle f(A)x,x\rangle\leq f(\langle Ax, x\rangle),
$$
\item[(ii)] there exists $\beta>0$ such that
$$
\frac{1}{2^{1-s}}\langle f(A)x,x\rangle-\beta\leq f(\langle Ax, x\rangle).
$$
\end{itemize}
\item [(4)] If $f$ is Godunova-Levin on $[m, M]$, then
\begin{itemize}
\item[(i)] there exists $\alpha>0$ such that
$$
\frac{1}{4\alpha}\langle f(A)x,x\rangle\leq f(\langle Ax, x\rangle),
$$
\item[(ii)] there exists $\beta>0$ such that
$$
\frac{1}{4}\langle f(A)x,x\rangle-\beta\leq f(\langle Ax, x\rangle).
$$
\end{itemize}
\item [(5)] If $f$ is $s$-Godunova-Levin on $[m, M]$, then
\begin{itemize}
\item[(i)] there exists $\alpha>0$ such that
$$
\frac{1}{2^{1+s}\alpha}\langle f(A)x,x\rangle\leq f(\langle Ax, x\rangle),
$$
\item[(ii)] there exists $\beta>0$ such that
$$
\frac{1}{2^{1+s}}\langle f(A)x,x\rangle-\beta\leq f(\langle Ax, x\rangle).
$$
\end{itemize}
\end{itemize}
\end{corollary}
\begin{proof}
Consider $h(t)=t$, $h(t)=1$, $h(t)=t^s$, $h(t)=\frac{1}{t}$, and $h(t)=\frac{1}{t^s}$ in parts (1)-(5), respectively
and note that
the function $\frac{h(t)}{t}$ is decreasing in each part and Corollary \ref{thebest}
implies the Jensen's coefficient is $2h(\frac{1}{2})$.
According Theorem \ref{cnrs-thm4}, we get the desired result in each part.
\end{proof}

\section{Applications}
In this section,
we obtain the Hermite-Hadamard's type inequality for $h$-convex functions.
Moreover, we obtain a multiple operator version of Theorem \ref{thm1} for $h$-convex functions.
In particular, one may reach a result for convex, $P$-class, $s$-convex, Godunova-Levin, and  $s$-Godunova-Levin functions.
% %%%%%%%%%%%%%%%%%%%%%%%%%%%%%%%%%%%%%%%%%%%%%%%%%%%%%%%%
%
\begin{corollary}\label{thm5}
Let the conditions of Theorem \ref{thm1} be satisfied and let $p$ and $q$
be non-negative numbers, with $p + q > 0$, for which
\begin{align*}
\langle Ax,x \rangle=\frac{pm+qM}{p+q}.
\end{align*}
Then,
\begin{align*}
\frac{1}{M_{(0,1)}(h)}f\Big(\frac{pm+qM}{p+q}\Big)\leq \langle f(A)x,x \rangle\leq M_{(0,1)}(h)\frac{pf(m)+qf(M)}{p+q}.
\end{align*}
\end{corollary}
%%%%%%%%%%%%%%%%%%%%%%%%%%%%%%%
\begin{proof}
By virtue of Theorem  \ref{thm1} and \ref{thm3} we reach
\begin{align*}%\label{eq1-thm5} %\frac{1}{2}
f\Big(\frac{pm+qM}{p+q}\Big)&=f(\langle Ax,x \rangle)\leq M_{(0,1)}(h)\langle f(A)x,x \rangle)\\
&\leq (M_{(0,1)}(h))^2\Big(\frac{M-\langle Ax,x\rangle}{M-m}f(m)+\frac{\langle Ax,x\rangle-m}{M-m}f(M)\Big)\\
&=(M_{(0,1)}(h))^2\frac{pf(m)+qf(M)}{p+q}.
\end{align*}
\end{proof}

We may consider a multiple operator version of Theorem \ref{thm1} as follows
and obtain some interesting corollaries.

% %%%%%%%%%%%%%%%%%%%%%%%%%%%%%%%%%%%%%%%%%%%%%%%%%%%%%%%%
%
\begin{theorem}\label{M-Mond-pclass}
Let $A_i$ be self-adjoint operators with $Sp(A_i)\subseteq [m,M]$ for some scalars
$m<M$ and $x_i\in \mathcal H$, $i\in\{1,...,n\}$ with $\sum_{i=1}^{n}||x_i||^2=1$.
If $f$ is $h$-convex on $[m,M]$, then
$$
f\Big(\sum_{i=1}^{n}\langle A_ix_i,x_i\rangle\Big)\leq M_{(0,1)}(h)\sum_{i=1}^{n}\langle f(A_i)x_i,x_i\rangle.
$$
\end{theorem}
\begin{proof}
We define
$$
\tilde{A}=\left(
\begin{array}{lll}
A_1 & \cdots & 0\\
\vdots & \ddots & \vdots\\
0 & \cdots & A_n
\end{array}
\right) \ \ \ \ \ \  \textrm{and} \ \ \ \ \ \
\tilde{x}=\left(
\begin{array}{l}
x_1\\
\vdots\\
x_n
\end{array}
\right).
$$
So, $Sp(\tilde{A})\subseteq [m,M]$, $||\tilde{x}||=1$, and
$$
f(\langle \tilde{A}\tilde{x},\tilde{x} \rangle)=f\Big(\sum_{i=1}^{n}\langle A_ix_i,x_i\rangle\Big),
$$
$$
\langle f(\tilde{A})\tilde{x}, \tilde{x}\rangle=\sum_{i=1}^{n}\langle f(A_i)x_i,x_i\rangle.
$$
In view of Theorem \ref{thm1} the result follows.
\end{proof}

% %%%%%%%%%%%%%%%%%%%%%%%%%%%%%%%%%%%%%%%%%%%%%%%%%%%%%%%%
%
\begin{corollary}\label{M-Mond-pclass-cor}
Let $A_i$ be self-adjoint operators with $Sp(A_i)\subseteq [m,M]$, $i\in\{1,...,n\}$ for some scalars
$m<M$. If $f$ is $h$-convex on $[m,M]$ and $p_i\geq 0$ with $\sum_{i=1}^{n}p_i=1$, then
$$
f\Big(\sum_{i=1}^{n}p_i\langle A_ix,x\rangle\Big)\leq M_{(0,1)}(h)\sum_{i=1}^{n}p_i\langle f(A_i)x,x\rangle
$$
for every $x\in \mathcal H$ with $||x||=1$.
\end{corollary}
\begin{proof}
By using Theorem \ref{M-Mond-pclass} and setting $x_i=\sqrt{p_i}x$, $i\in\{1,...,n\}$ we can reach the result.
\end{proof}

%%%%%%%%%%%%%%%%%%%%%%%%%%%%%%%
%The following corollary is also of interest.
%In the following theorem, we make a refinement of Corollary \ref{M-Mond-pclass-cor}.

%%%%%%%%%%%%%%%%%%%%%%%%%%%%%%%
\begin{corollary} %\label{cor4}
Let the conditions of Theorem \ref{M-Mond-pclass} be satisfied.
\begin{itemize}
\item[(i)] There exists $\alpha>0$ such that
\begin{equation}
\frac{1}{M_{(0,1)}(h)\alpha}\sum_{i=1}^{n}\langle f(A_i)x_i,x_i\rangle\leq f(\langle \sum_{i=1}^{n}A_ix_i, x_i\rangle).
\end{equation}
\item[(ii)] There exists $\beta>0$ such that
\begin{equation}%\label{eq0-cor4}
\frac{1}{M_{(0,1)}(h)}\sum_{i=1}^{n}\langle f(A_i)x_i,x_i\rangle -\beta\leq f(\langle \sum_{i=1}^{n}A_ix_i, x_i\rangle).
\end{equation}
\end{itemize}
\end{corollary}
%%%%%%%%%%%%%%%%%%%%%%%%%%%%%%%
%%%%%%%%%%%%%%%%%%%%%%%%%%%%%%%
\begin{proof}
Consider $\tilde{A}$ and $\tilde{x}$ as in the proof of Theorem \ref{M-Mond-pclass} and
Apply Theorem \ref{cnrs-thm4}.
\end{proof}

%{\bf Acknowledgement:}
%The first author was supported by Payame Noor University under grant number~46384.

\section{Declarations}
We remark that the potential conflicts of interest and
data sharing not applicable to this article and no data sets were generated during the current study.

%\begin{acknowledgements}
%The first author was supported by Payame Noor University under grant number~46384.
%\end{acknowledgements}

% BibTeX users please use one of
%\bibliographystyle{spbasic}      % basic style, author-year citations
%\bibliographystyle{spmpsci}      % mathematics and physical sciences
%\bibliographystyle{spphys}       % APS-like style for physics
%\bibliography{}   % name your BibTeX data base

% Non-BibTeX users please use

\end{document}